\documentclass[11pt]{amsart}
\usepackage[margin=1.1in]{geometry}
\usepackage[dvipsnames]{xcolor}
\usepackage{amsmath, amsfonts, amssymb, amsthm}
\usepackage[unicode]{hyperref}
\usepackage[all]{hypcap}
\usepackage{fancyhdr}
\usepackage{bbm}
\usepackage{thm-restate}
\usepackage[mathlines]{lineno}

\modulolinenumbers[1]

\theoremstyle{plain}
\newtheorem{theorem}{Theorem}[section]

\newtheorem{lemma}[theorem]{Lemma}
\newtheorem{proposition}[theorem]{Proposition}
\newtheorem{corollary}[theorem]{Corollary}

\theoremstyle{definition}
\newtheorem{definition}{Definition}[section]

\newtheorem*{remark}{Remark}
\newtheorem*{example}{Example}

\newtheorem{theirtheorem}{Theorem}

\newcommand{\restatementtitle}{}
\newcommand{\restatetheorem}[1]{%
    \begingroup
    \renewcommand{\restatementtitle}{Restated}%
    #1*%
    \endgroup
}

\newcommand{\tensor}{\otimes}
\newcommand{\Cay}{\text{Cay}}
\newcommand{\Mat}{\text{Mat}}
\newcommand{\R}{\mathbb{R}}
\newcommand{\C}{\mathbb{C}}
\newcommand{\Z}{\mathbb{Z}}

\newcommand{\oA}{{\overline A}}
\newcommand{\oB}{{\overline B}}
\newcommand{\oC}{{\overline C}}

\newcommand{\cP}{\mathcal{P}}
\newcommand{\cQ}{\mathcal{Q}}

\newcommand{\oc}{{\overline c}}
\newcommand{\og}{{\overline g}}

\newcommand{\ou}{{\overline u}}

\newcommand{\seq}{\subseteq}
\newcommand{\stm}{\setminus}

\newcommand{\pgh}[1][G/H]{\phi_{#1}}

\title[Q-Polynomial Graphs and Adjacency Algebras]{\textbf{On $Q$-Polynomial Graphs and Entrywise Closure of Their Adjacency Algebras}}
\author{Ivy Liu}
\address{Department of Mathematics, University of Wisconsin--Madison, Madison, Wisconsin, USA}
\email{ivy.liu@wisc.edu}
\date{}

\begin{document}

\begin{abstract}
We study the $Q$-polynomial property for general graphs introduced
by P. Terwilliger. We prove that the Cartesian product of two
$Q$-polynomial graphs is $Q$-polynomial whenever their eigenvalues,
in suitable $Q$-polynomial orderings, form arithmetic progressions
with the same common difference. We use this result to construct an infinite family of
$Q$-polynomial graphs whose adjacency algebras are not
closed under entrywise multiplication, disproving a conjecture
of P. Terwilliger. We then propose a refined conjecture:
every $Q$-polynomial graph whose adjacency algebra is not closed
under entrywise multiplication admits two nontrivial equitable
partitions such that each cell of one partition meets each cell of the other
in exactly one vertex. We prove this refined conjecture for
abelian Cayley graphs.
\end{abstract}

\maketitle


\section{Introduction}

For distance-regular graphs, the $Q$-polynomial property has been studied for several decades, starting with Delsarte's work on association schemes~\cite{Delsarte1973}. Later, Bannai and Ito further developed this theory in their monograph~\cite{BannaiIto1984}. More recently, the $Q$-polynomial property for general graphs that need not be distance-regular has been introduced by Terwilliger~\cite{Terwilliger2023DistanceRegular}. In this paper, we consider the $Q$-polynomial property in this general setting. The precise definition for $Q$-polynomial is given in Section~\ref{subsec:q_poly}.

Our first main result is a method to construct $Q$-polynomial graphs using Cartesian products. 

\begin{restatable}[name=\restatementtitle]{theorem}{cartesianQPtheorem}\label{thm:cartesian_QP}
    Suppose $\Gamma_1$ and $\Gamma_2$ are $Q$-polynomial graphs. Let
    $\lambda_0^{(1)}, \ldots, \lambda_s^{(1)}$ and $\lambda_0^{(2)},\ldots,\lambda_t^{(2)}$
    be $Q$-polynomial orderings of their eigenvalues. Assume that these eigenvalue sequences are in arithmetic progression with the same common difference. Then the Cartesian product $\Gamma_1 \times \Gamma_2$ is $Q$-polynomial. 
    
    More precisely, for $0\le k\le s+t$, define 
    \[
        \theta_k = \lambda^{(1)}_i + \lambda^{(2)}_j
        \quad\text{where } k = i+j, \quad 0 \leq i \leq s, \quad 0 \leq j \leq t.
    \]
    Then the sequence $\theta_0,\theta_1,\ldots,\theta_{s+t}$ is a $Q$-polynomial ordering of the eigenvalues of $\Gamma_1 \times \Gamma_2$.
\end{restatable}

We illustrate Theorem~\ref{thm:cartesian_QP} with an example. To describe the example, we bring in some notation. For a graph $\Gamma$
and an integer $n\geq 1$ define the Cartesian product 
\[
        \Gamma^{n}
        :=
        \underbrace{\Gamma\times\cdots\times\Gamma}_{n\text{ copies}}.
    \]
Fix an integer $r\geq 3$. Our example involves the complete graph $K_r$ and the complete bipartite graph $K_{r,r}$. In Section~\ref{subsec:q_poly}, we show that $r-1,-1$ is a $Q$-polynomial ordering of the eigenvalues of $K_r$ and that $r,0,-r$ is a $Q$-polynomial ordering of the eigenvalues of $K_{r,r}$. These two sequences are arithmetic progressions with the same common difference $-r$. By repeated application of Theorem~\ref{thm:cartesian_QP}, the Cartesian
product
\begin{equation}
    K_r^m \times K_{r,r}^n
    \tag{*}\label{eq:product-example}
\end{equation}
is $Q$-polynomial for all integers $m,n\geq 1$.

We will use the example \eqref{eq:product-example} to disprove a conjecture proposed by P. Terwilliger~\cite[Conjecture~20.4]{Terwilliger2023DistanceRegular}. The conjecture states that the adjacency algebra of every $Q$-polynomial graph is closed under entrywise
multiplication. In Theorem~\ref{thm:counterexample} we will show that every graph of the form
\eqref{eq:product-example} is a counterexample to \cite[Conjecture~20.4]{Terwilliger2023DistanceRegular}.

In Theorem~\ref{thm:cartesian_QP} we gave a method to construct $Q$-polynomial graphs using Cartesian products. Later in this paper, we will give a similar method to construct $Q$-polynomial graphs using tensor products. The resulting $Q$-polynomial graphs will potentially yield further counterexamples to \cite[Conjecture~20.4]{Terwilliger2023DistanceRegular}. Moreover, if $\Gamma$ is any $Q$-polynomial graph that is a counterexample to \cite[Conjecture~20.4]{Terwilliger2023DistanceRegular} and the complement of $\Gamma$ is connected, then the complement is again a counterexample, since $\Gamma$ and its complement have the same eigenspaces and the same adjacency algebra. Motivated by these constructions, we give an adjusted version of \cite[Conjecture~20.4]{Terwilliger2023DistanceRegular} for which there are no counterexamples known to us.


Recall that a partition \(\mathcal P=\{P_1,\ldots,P_t\}\) of the vertex set of a graph is \emph{equitable} if, for every \(i,j\), the number of neighbors in \(P_j\) of a vertex \(x\in P_i\) is independent of the choice of \(x\). A partition is called \emph{trivial} if it has a single part. Observe that an equitable partition of a graph is also an equitable partition of its complement.

The counterexamples constructed using Cartesian products and the potential counterexamples suggested by tensor products share a common property: their vertex sets admit two nontrivial equitable partitions such that each cell of one partition meets each cell of the other in exactly one vertex. To describe these partitions, write the vertex set of a product as $X_1\times X_2$. Fixing $x\in X_1$ gives a cell $\{x\}\times X_2$ of the first partition; fixing $y\in X_2$ gives a cell $X_1\times\{y\}$ of the second. These two cells meet in the single vertex $(x,y)$. The same partitions remain equitable in the complement, so the complements share this property as well. This observation motivates the following conjecture.

\begin{restatable}[name=\restatementtitle]{conjecture}{refinedconjecture}\label{conj:new:intro}\label{conj:new}
    If the adjacency algebra of a $Q$-polynomial graph is not closed under entrywise multiplication, then its vertex set admits two nontrivial equitable partitions such that each cell of one partition meets each cell of the other in exactly one vertex.
\end{restatable}

For abelian Cayley graphs, character theory provides an additive formulation of the $Q$-polynomial property: the eigenvectors of the abelian Cayley graph $\Cay(G,S)$ can be indexed by elements of $G$, and entrywise multiplication of eigenvectors corresponds to addition in $G$. When the group $G$ is decomposable, its decomposition yields the required equitable partitions. When $G$ is indecomposable, this additive formulation allows us to apply Kemperman's structure theorem from additive combinatorics to show that the adjacency algebra coincides with that of a distance-regular graph on the same vertex set and is therefore closed under entrywise multiplication.

\begin{restatable}[name=\restatementtitle]{theorem}{abeliancayleytheorem}\label{thm:abelian_cayley}
    Conjecture~\ref{conj:new:intro} holds true for abelian Cayley graphs.
\end{restatable}

Section~\ref{sec:preliminaries} reviews the necessary preliminaries. Section~\ref{sec:counterexample} presents a method for constructing $Q$-polynomial graphs (Theorem~\ref{thm:cartesian_QP}), constructs a family of counterexamples (Theorem~\ref{thm:counterexample}) to Conjecture~\cite[Conjecture~20.4]{Terwilliger2023DistanceRegular}, and formulates the refined conjecture (Conjecture~\ref{conj:new}). Section~\ref{sec:new_conj} establishes Conjecture~\ref{conj:new} for abelian Cayley graphs (Theorem~\ref{thm:abelian_cayley}).

\section{Preliminaries}\label{sec:preliminaries}

Throughout this paper, $\Gamma=(X,\mathcal R)$ denotes a finite, connected, undirected graph without multiple edges, with nonempty vertex set $X$.

\subsection{Entrywise multiplication}\label{subsec:entrywise_mul}
Recall the nonempty finite set $X$. Let $\Mat_X(\R)$ denote the algebra of real
matrices whose rows and columns are indexed by $X$, and let $V:=\R^X$ denote
the space of real column vectors indexed by $X$.
The algebra $\Mat_X(\R)$ acts on $V$ by matrix multiplication.

We now introduce entrywise multiplication. For
$B,C\in\Mat_X(\R)$, define $B\circ C\in\Mat_X(\R)$ by
\[
    (B\circ C)_{x,y}=B_{x,y}C_{x,y}
    \qquad (x,y\in X).
\]
The operation $\circ$ is called the \emph{entrywise multiplication}.

\begin{proposition} \label{prop:tensor_circ}
    For matrix $A,B,C,D$ where $A,C$ has same dimension and $B,D$ has same dimension,
    \[(A \tensor B) \circ (C \tensor D) = (A \circ C) \tensor (B \circ D)\]
\end{proposition}

Similarly, for $u,v\in V$, define $u \circ v \in V$ as
\[
    (u\circ v)_x=u_xv_x
    \qquad (x\in X).
\]
For subspaces $U,W\subseteq V$, define the subspace $U \circ W$ of $V$ to be 
\[
    U\circ W
    =
    \operatorname{span}_{\R}
    \{u\circ w:u\in U,\ w\in W\}.
\]
For an integer $\ell\geq 0$, define
\[
    U^{\circ\ell}
    =
    \underbrace{U\circ U\circ\cdots\circ U}_{\ell\text{ copies}}.
    \qquad
\]
For $\ell = 0$, we interpret $U^{\circ 0}$ as $\operatorname{span}_\R\{\mathbf 1\}$, where $\mathbf 1$ denotes the all-ones vector.

\subsection{General graphs}\label{sec:general graph}
Recall the graph $\Gamma = (X, \mathcal R)$ and let $A\in\Mat_X(\R)$ denote the $0/1$ adjacency matrix of $\Gamma$. We call the algebra generated by $A$ the adjacency algebra of $\Gamma$ and denote it by $M$. The space $V$, defined in Section \ref{subsec:entrywise_mul}, is the standard $M$-module then, with $M$ acting by matrix multiplication.
Let $\Lambda$ be the set of distinct eigenvalues of $A$, and put $d=|\Lambda|-1$.
For each $\lambda\in\Lambda$, let $V_\lambda$ be the corresponding eigenspace of $A$.
Since $M$ consists of mutually commuting real symmetric matrices, these are
its common eigenspaces, and $V=\bigoplus_{\lambda\in\Lambda}V_\lambda$.
Define $E_\lambda\in\Mat_X(\R)$ to be the orthogonal projector
onto $V_\lambda$, so that $E_\lambda V=V_\lambda$.
The matrices $\{E_\lambda:\lambda\in\Lambda\}$ form a basis of $M$ and are called its
\emph{primitive idempotents} with respect to ordinary matrix multiplication.
In particular,
\[
    A=\sum_{\lambda\in\Lambda}\lambda E_\lambda.
\]

Define an equivalence relation on
$X\times X$ by
\[
    (x,y)\sim(x',y')
    \quad\Longleftrightarrow\quad
    B_{x,y}=B_{x',y'}
    \text{ for every }B\in M.
\]
In particular, since $\{E_\lambda:\lambda\in\Lambda\}$ forms a basis of $M$, this is equivalent to
$(E_\lambda)_{x,y}=(E_\lambda)_{x',y'}$ for every $\lambda\in\Lambda$.
Let $R_0,\ldots,R_\delta$ be the resulting equivalence classes. For
$0\leq i\leq \delta$, define the characteristic matrix $A_i$ of $R_i$ by
\[
    (A_i)_{x,y}
    =
    \begin{cases}
        1, & (x,y)\in R_i,\\
        0, & (x,y)\notin R_i.
    \end{cases}
\]

\begin{proposition}\label{prop:@Ai=@Ei}
    $M$ is closed under $\circ$ if and only if $d=\delta$.
\end{proposition}
\begin{proof}
    This follows from \cite[Theorem~2.6.1(i)]{BCN1989}.
\end{proof}

When these equivalent conditions hold, the matrices $A_0,\ldots,A_\delta$
form a basis of $M$ consisting of pairwise orthogonal idempotents with
respect to $\circ$. They are therefore the primitive idempotents of $M$
with respect to entrywise multiplication.

\subsection{Distance-regular graphs}
Let $\partial(x,y)$ denote the distance of two vertices $x,y \in X$, and let $D$ denote the diameter of
$\Gamma$, that is, the maximum distance between any two vertices of $\Gamma$.
For $x\in X$ and $0\leq i\leq D$, define
\[
    \Gamma_i(x)=\{y\in X:\partial(x,y)=i\}.
\]

We now recall the definition of a distance-regular graph.

\begin{definition}[Distance-regular graph]
    The graph $\Gamma$ is called \emph{distance-regular} if, for all
    $0\leq h,i,j\leq D$ and all vertices $x,y\in X$ with
    $\partial(x,y)=h$, the number
    \[
        p_{ij}^h
    =
    \lvert \Gamma_i(x)\cap\Gamma_j(y)\rvert
    \]
    depends only on $h,i,j$, and not on the choice of $x,y$. The numbers
    $p_{ij}^h$ are called the \emph{intersection numbers} of $\Gamma$.
\end{definition}
\begin{example}[Complete graph]
    The complete graph $K_r$ has $r$ vertices, with every pair of distinct vertices
    adjacent. It is distance-regular with
    \[
        p_{00}^0=1,\qquad p_{11}^0=r-1,\qquad
        p_{01}^1=p_{10}^1=1,\qquad p_{11}^1=r-2,
    \]
    and all other intersection numbers equal to zero.
\end{example}

\begin{example}[Strongly regular graph]
    $\Gamma$ is called \emph{strongly regular} if it is distance-regular with diameter $2$.
\end{example}

\begin{example}[Cycle]
    For $n\geq 3$, the cycle $C_n$ is the graph with vertex set $\Z_n$ of integers
    modulo $n$, in which two vertices are adjacent exactly when their difference
    is $\pm1$ modulo $n$. Let $D=\lfloor n/2\rfloor$.
    For $0\leq i\leq D$, set
    \[
        S_i=\{\pm i\},
    \]
    where repeated elements are counted only once. The cycle $C_n$ is
    distance-regular with
    \[
        p_{ij}^h=|S_i\cap(h+S_j)|.
    \]
\end{example}

When $\Gamma$ is distance-regular, $M$ is also called its \emph{Bose--Mesner algebra}.
For $0\leq i\leq D$, define the $i^{\text{th}}$ distance matrix $A_i\in\Mat_X(\R)$ by
\[
    (A_i)_{x,y}
    =
    \begin{cases}
        1, & \partial(x,y)=i,\\
        0, & \partial(x,y)\neq i,
    \end{cases}
    \qquad x,y\in X.
\]
Here $A_1=A$. By \cite[Section~2 and Corollary~3.4(i)]{Terwilliger2023DistanceRegular},
the matrices $A_0,\ldots,A_D$ form a basis of $M$.
In particular, this definition agrees with the definition in
Section~\ref{sec:general graph} when the classes are indexed by distance.

These distance matrices satisfy
\[
    A_i\circ A_j=\delta_{ij}A_i
    \qquad (0\leq i,j\leq D).
\]
It follows that the Bose--Mesner algebra of every distance-regular
graph is closed under entrywise multiplication. In particular, this
holds for complete graphs, connected strongly regular graphs, and
cycles.

\subsection{$Q$-polynomial Property} \label{subsec:q_poly}
We now recall the definition of the $Q$-polynomial property for a
general graph.
\begin{definition}[$Q$-polynomial graph]
Let
$\lambda_0,\lambda_1,\ldots,\lambda_d$
be an ordering of $\Lambda$, and let $E_i:=E_{\lambda_i}$ and
$V_i:=V_{\lambda_i}$ be the corresponding primitive idempotents and eigenspaces.
The ordering $V_0,V_1,\ldots,V_d$ is called \emph{$Q$-polynomial} if
\[
    \sum_{\ell=0}^i V_\ell
    =
    \sum_{\ell=0}^i V_1^{\circ\ell}
    \qquad (0\leq i\leq d).
\]
In particular, $V_0 = \operatorname{span}\{\mathbf 1\}$. The graph $\Gamma$ is called \emph{$Q$-polynomial} if there exists at least one $Q$-polynomial ordering of the eigenspaces of $A$.
\end{definition}

This definition is taken from
\cite[Definitions~20.1 and~20.2]{Terwilliger2023DistanceRegular}.
For distance-regular graphs, it agrees with the usual
$Q$-polynomial property by
\cite[Theorem~12.1]{Terwilliger2023DistanceRegular}.

\begin{example}
    Let $r\geq 2$. The complete graph $K_r$ is trivially $Q$-polynomial since it only has two distinct eigenvalues.
\end{example}

\begin{example}[Strongly regular graphs]
    Connected and non-complete strongly regular graphs are precisely the distance-regular graphs with diameter $2$, and thus are $Q$-polynomial; see
    \cite[Page~5]{DeBruynVanhove2013Inequalities}. In particular, the complete bipartite graph \(K_{r,r}\) has \(Q\)-polynomial eigenvalue ordering \(r,0,-r\) for \(r\ge 2\). Indeed, the \((-r)\)-eigenspace is spanned by \(w=(\mathbf 1_r,-\mathbf 1_r)\), and \(w\circ w=\mathbf 1\) rules out the ordering \(r,-r,0\).
\end{example}

\subsection{Cartesian Product of Graphs}\label{subsec:cartesian_product}

\begin{definition}[Cartesian product of graphs]
    Let $\Gamma_1 = (X_1, \mathcal R_1)$ and $\Gamma_2 = (X_2, \mathcal R_2)$ be graphs. The \emph{Cartesian product} of
    $\Gamma_1$ and $\Gamma_2$, denoted by $\Gamma_1 \times \Gamma_2$, is the
    graph with vertex set $X_1 \times X_2$, where two vertices $(x_1,x_2)$ and $(y_1,y_2)$ are adjacent if and only if
    either
    \[
        x_1=y_1
        \quad\text{and}\quad
        x_2\sim_{\Gamma_2}y_2,
    \]
    or
    \[
        x_2=y_2
        \quad\text{and}\quad
        x_1\sim_{\Gamma_1}y_1.
    \]
\end{definition}

Assume until further notice that $\Gamma = \Gamma_1 \times \Gamma_2$. Let $A_1$ and $A_2$ be the adjacency matrices of
$\Gamma_1$ and $\Gamma_2$, respectively. Under the natural identification
\[
    \R^X
    \cong
    \R^{X_1}\otimes\R^{X_2},
\]
the adjacency matrix of $\Gamma_1\times\Gamma_2$ is
\[
    A
    =
    A_1\otimes I_{|X_2|}
    +
    I_{|X_1|}\otimes A_2.
\]

Let $\Lambda^{(1)}$ and $\Lambda^{(2)}$ be the sets of distinct eigenvalues of
$\Gamma_1$ and $\Gamma_2$, respectively.
For $h\in\{1,2\}$ and $\lambda\in\Lambda^{(h)}$, let
$U_\lambda^{(h)}$ be the corresponding eigenspace of $A_h$, and let
$E_\lambda^{(h)}$ be the orthogonal projector onto $U_\lambda^{(h)}$.
Thus,
\[
    A_1
    =
    \sum_{\lambda\in\Lambda^{(1)}}\lambda E_\lambda^{(1)},
    \qquad
    A_2
    =
    \sum_{\mu\in\Lambda^{(2)}}\mu E_\mu^{(2)}.
\]
It follows that
\[
    A
    =
    \sum_{\lambda\in\Lambda^{(1)}}\sum_{\mu\in\Lambda^{(2)}}
    (\lambda+\mu)(E_\lambda^{(1)}\otimes E_\mu^{(2)}).
\]

For any $u\in U_\lambda^{(1)}$ and $w\in U_\mu^{(2)}$, we have
\[
    A(u\otimes w)
    =
    (\lambda+\mu)(u\otimes w).
\]
Consequently, the set of eigenvalues of $\Gamma$ is
\[
    \Lambda=\{\lambda+\mu:\lambda\in\Lambda^{(1)},\ \mu\in\Lambda^{(2)}\}.
\]

For each $\theta\in\Lambda$, the corresponding eigenspace is
\[
    V_\theta
    =
    \bigoplus_{\substack{\lambda\in\Lambda^{(1)},\ \mu\in\Lambda^{(2)}\\
                        \lambda+\mu=\theta}}
    U_\lambda^{(1)}\otimes U_\mu^{(2)},
\]
and the corresponding primitive idempotent is
\[
    E_\theta
    =
    \sum_{\substack{\lambda\in\Lambda^{(1)},\ \mu\in\Lambda^{(2)}\\
                    \lambda+\mu=\theta}}
    E_\lambda^{(1)}\otimes E_\mu^{(2)}.
\]
Therefore,
\[
    A
    =
    \sum_{\theta\in\Lambda}\theta E_\theta.
\]
Moreover, the multiplicity of $\theta$ is
\[
    \dim V_\theta
    =
    \sum_{\substack{\lambda\in\Lambda^{(1)},\ \mu\in\Lambda^{(2)}\\
                    \lambda+\mu=\theta}}
    \dim U_\lambda^{(1)}\dim U_\mu^{(2)}.
\]
In particular, if all the sums $\lambda+\mu$ are distinct, then the
primitive idempotents of $\Gamma_1\times\Gamma_2$ are precisely $E_\lambda^{(1)}\otimes E_\mu^{(2)}$ where $\lambda\in\Lambda^{(1)}$ and $\mu\in\Lambda^{(2)}$.

Now let $\Gamma$ be an arbitrary graph. Recall that for a positive integer $n$, its $n$-fold Cartesian power is
\[
    \Gamma^{n}
    :=
    \underbrace{\Gamma\times\cdots\times\Gamma}_{n\text{ copies}}.
\]
The adjacency matrix of $\Gamma^{n}$ is
\[
    A_{\Gamma^{n}}
    =
    \sum_{t=1}^{n}
    I^{\otimes(t-1)}
    \otimes A
    \otimes I^{\otimes(n-t)},
\]
where $I=I_{|X|}$. The set of eigenvalues of $\Gamma^{n}$ is $\{\lambda^{(1)}+\cdots+\lambda^{(n)}:\lambda^{(1)},\ldots,\lambda^{(n)}\in\Lambda\}$. Here the superscripts on $\lambda^{(1)},\ldots,\lambda^{(n)}$ label the tensor factors; repetitions are allowed. The
primitive idempotent corresponding to an eigenvalue $\theta$ of
$\Gamma^{n}$ is
\[
    \sum_{\substack{\lambda^{(1)},\ldots,\lambda^{(n)}\in\Lambda\\
                    \lambda^{(1)}+\cdots+\lambda^{(n)}=\theta}}
    E_{\lambda^{(1)}}\otimes\cdots\otimes E_{\lambda^{(n)}}.
\]

\section{Constructing $Q$-polynomial graphs via Cartesian products}\label{sec:counterexample}
In this section, we provide a new method to construct $Q$-polynomial graphs through Cartesian product, and use this method to construct a family of counterexamples to \cite[Conjecture~20.4]{Terwilliger2023DistanceRegular}, i.e., a family of Q-polynomial graphs whose adjacency algebra is not closed under $\circ$.

The following theorem shows that the Cartesian product of two $Q$-polynomial graphs, with eigenvalues forming arithmetic progressions with the same gap in their respective $Q$-polynomial orderings, is a $Q$-polynomial graph.

\restatetheorem{\cartesianQPtheorem}

\begin{proof}
    Since the two eigenvalue sequences are arithmetic progressions with the same common difference,
    \[
        \lambda_i^{(1)}+\lambda_j^{(2)}
        =\lambda_{i'}^{(1)}+\lambda_{j'}^{(2)}
        \quad\Longleftrightarrow\quad i+j=i'+j'.
    \]
    Thus the eigenvalues of $\Gamma_1\times\Gamma_2$ can be ordered as
    $\theta_0,\ldots,\theta_{s+t}$, where $\theta_k=\lambda_i^{(1)}+\lambda_j^{(2)}$ for $i+j=k$.
    Let $V_k:=V_{\theta_k}$ and $E_k:=E_{\theta_k}$ be the corresponding eigenspaces and primitive idempotents.
    For $h\in\{1,2\}$, let $U_i^{(h)}$ denote the eigenspace of $\Gamma_h$ corresponding to $\lambda_i^{(h)}$.
    Set $U_i^{(1)}=0$ for $i>s$ and $U_j^{(2)}=0$ for $j>t$.
    By Section~\ref{subsec:cartesian_product},
    \[V_k = \bigoplus_{i+j=k}U_i^{(1)} \tensor U_j^{(2)}.\]
    We shall prove that this ordering is $Q$-polynomial.

    By a direct calculation and applying Proposition \ref{prop:tensor_circ}, one can show that
    \[V_1^{\circ k} = \sum_{i+j=k}\bigl(U_1^{(1)}\bigr)^{\circ i} \tensor \bigl(U_1^{(2)}\bigr)^{\circ j},\]
    \[\sum_{l=0}^k V_1^{\circ l} = \sum_{i+j \leq k} \bigl(U_1^{(1)}\bigr)^{\circ i} \tensor \bigl(U_1^{(2)}\bigr)^{\circ j}.\]
    Moreover,
    \[\sum_{l=0}^k V_l = \sum_{i+j \leq k} U_i^{(1)} \tensor U_j^{(2)}.\]
    In each factor, both cumulative sums in the $Q$-polynomial condition equal the full space at the last index of the ordering, so they remain equal at all larger indices.
    Thus, for $i+j \leq k$, we have
    \[
    \begin{aligned}
        \bigl(U_1^{(1)}\bigr)^{\circ i}\tensor\bigl(U_1^{(2)}\bigr)^{\circ j}
        &\subseteq \left(\sum_{p=0}^i\bigl(U_1^{(1)}\bigr)^{\circ p}\right)
        \tensor\left(\sum_{q=0}^j\bigl(U_1^{(2)}\bigr)^{\circ q}\right)\\
        &=\left(\sum_{p=0}^i U_p^{(1)}\right)\tensor\left(\sum_{q=0}^j U_q^{(2)}\right)
        \subseteq\sum_{l=0}^k V_l,
    \end{aligned}
    \]
    \[
    \begin{aligned}
        U_i^{(1)}\tensor U_j^{(2)}
        &\subseteq\left(\sum_{p=0}^i U_p^{(1)}\right)\tensor\left(\sum_{q=0}^j U_q^{(2)}\right)\\
        &=\left(\sum_{p=0}^i\bigl(U_1^{(1)}\bigr)^{\circ p}\right)
        \tensor\left(\sum_{q=0}^j\bigl(U_1^{(2)}\bigr)^{\circ q}\right)
        \subseteq\sum_{l=0}^k V_1^{\circ l}.
    \end{aligned}
    \]
    Therefore, $\sum_{l=0}^k V_1^{\circ l} = \sum_{l=0}^k V_l$ and thus $\{V_k\}$ is $Q$-polynomial.

\end{proof}

Theorem~\ref{thm:cartesian_QP} has an analogue in which Cartesian products and arithmetic progressions are replaced by tensor products and geometric progressions, respectively. We first define the tensor product of graphs.

\begin{definition}[Tensor product of graphs]\label{def:tensor_product}
    Let $\Gamma_1 = (X_1, \mathcal R_1)$ and $\Gamma_2 = (X_2, \mathcal R_2)$ be graphs. The \emph{tensor product} of
    $\Gamma_1$ and $\Gamma_2$, denoted by $\Gamma_1 \tensor \Gamma_2$, is the
    graph with vertex set $X_1 \times X_2$, where two vertices $(x_1,x_2)$ and $(y_1,y_2)$ are adjacent if and only if
    \[
        x_1\sim_{\Gamma_1}y_1
        \quad\text{and}\quad
        x_2\sim_{\Gamma_2}y_2.
    \]
\end{definition}

\begin{theorem}\label{thm:tensor_QP}
    Suppose $\Gamma_1$ and $\Gamma_2$ are two nontrivial $Q$-polynomial graphs whose tensor product is connected. Let
    $\lambda_0^{(1)}, \ldots, \lambda_s^{(1)}$ and $\lambda_0^{(2)},\ldots,\lambda_t^{(2)}$
    be $Q$-polynomial orderings of their eigenvalues. Assume that these eigenvalue sequences are in geometric progression with the same common ratio. Then the tensor product $\Gamma_1 \tensor \Gamma_2$ is $Q$-polynomial.

    More precisely, for $0\le k\le s+t$, define
    $$\theta_k = \lambda^{(1)}_i\lambda^{(2)}_j$$
    where
    \[k = i+j, \quad 0 \leq i \leq s, \quad 0 \leq j \leq t.\]
    Then the sequence $\theta_0,\theta_1,\ldots,\theta_{s+t}$ is a $Q$-polynomial ordering of the eigenvalues of $\Gamma_1 \tensor \Gamma_2$.
\end{theorem}
\begin{proof}
    Let $\rho$ be the common ratio, and let $U_i^{(h)}$ be the eigenspace of $\Gamma_h$ for $\lambda_i^{(h)}$. If $A_1$ and $A_2$ are the adjacency matrices of the factors, then the adjacency matrix of $\Gamma_1\tensor\Gamma_2$ is $A_1\tensor A_2$. It acts on $U_i^{(1)}\tensor U_j^{(2)}$ as multiplication by
    \[
        \lambda_i^{(1)}\lambda_j^{(2)}
        =\lambda_0^{(1)}\lambda_0^{(2)}\rho^{i+j}.
    \]
    The factors are nontrivial regular graphs, so their zero traces and spectral bounds give $-1\leq\rho<0$. If $\rho=-1$, then $U_0^{(1)}\tensor U_0^{(2)}$ and $U_1^{(1)}\tensor U_1^{(2)}$ both lie in the eigenspace for the degree of the product, contradicting its connectedness. Thus $-1<\rho<0$, so $\theta_0,\ldots,\theta_{s+t}$ are distinct, with corresponding eigenspaces
    \[
        V_k=\bigoplus_{i+j=k}U_i^{(1)}\tensor U_j^{(2)}.
    \]
    By a calculation similar to that in the proof of Theorem~\ref{thm:cartesian_QP}, one can show that
    \[
        \sum_{\ell=0}^k V_1^{\circ\ell}
        =\sum_{i+j\leq k}U_i^{(1)}\tensor U_j^{(2)}
        =\sum_{\ell=0}^k V_\ell
        \qquad (0\leq k\leq s+t).
    \]
    Hence this ordering is $Q$-polynomial.
\end{proof}


We now return to Cartesian products and use Theorem~\ref{thm:cartesian_QP} to construct a family of counterexamples to \cite[Conjecture~20.4]{Terwilliger2023DistanceRegular}. We begin with the following corollary.

\begin{corollary}\label{cor:is_QP}
    Let $r\geq 2$ be an integer and let $m,n\geq 0$ be integers but not both zero. Then $K_r^{m} \times K_{r,r}^{n}$ is $Q$-polynomial.
    A factor with exponent zero is omitted.
\end{corollary}
\begin{proof}
    Set $\Gamma_1=K_r$ and $\Gamma_2=K_{r,r}$. Choose the $Q$-polynomial eigenvalue orderings
    \[\lambda_0^{(1)} = r-1, \lambda_1^{(1)} = -1\]
    and
    \[\lambda_0^{(2)} = r, \lambda_1^{(2)} = 0, \lambda_2^{(2)} = -r\]
    for $\Gamma_1$ and $\Gamma_2$, respectively. Both orderings have common difference $-r$, and each application of Theorem~\ref{thm:cartesian_QP} preserves this common difference. The result follows by induction on $m+n$.
\end{proof}

We now show that $K_r^m \times K_{r,r}^n$ gives a counterexample to \cite[Conjecture~20.4]{Terwilliger2023DistanceRegular} for $r\geq 3$, $m\geq 0$, $n\geq 1$, and $m+n\geq 2$. For a $Q$-polynomial graph $\Gamma=(X,\mathcal R)$ with the $Q$-polynomial ordering of primitive idempotents $E_0,\ldots,E_d$, define
\[
    N(\Gamma):=\#\{(E_1)_{x,y}:x,y\in X\}.
\]
This notation always refers to the chosen ordering.

\begin{proposition}\label{prop:N_cartesian}
    Under the hypotheses of Theorem~\ref{thm:cartesian_QP}, with $s,t\geq 1$ and the induced $Q$-polynomial ordering on the product, we have
    \[
        N(\Gamma_1\times\Gamma_2)\geq N(\Gamma_1)+N(\Gamma_2)-1.
    \]
\end{proposition}
\begin{proof}
    Let $X_1,X_2$ be the vertex sets of the factors. The second primitive idempotent of the product is
    \[
        E_1=E_1^{(1)}\tensor E_0^{(2)}+E_0^{(1)}\tensor E_1^{(2)}.
    \]
    Since every entry of $E_0^{(h)}$ is $1/|X_h|$, the entries of $|X_1||X_2|E_1$ are precisely the sums $a_i+b_j$, where
    $a_1<\cdots<a_p$ and $b_1<\cdots<b_q$ are the distinct entries of $|X_1|E_1^{(1)}$ and $|X_2|E_1^{(2)}$, respectively. The chain
    \[
        a_1+b_1<\cdots<a_p+b_1<a_p+b_2<\cdots<a_p+b_q
    \]
    gives at least $p+q-1=N(\Gamma_1)+N(\Gamma_2)-1$ distinct entries.
\end{proof}

For integers $r\geq 3$ and $m,n\geq 0$ with $m+n\geq 1$, put
\[
    N(m,n,r):=N(K_r^m\times K_{r,r}^n),
\]
using the ordering from Corollary~\ref{cor:is_QP}, and let $d(m,n,r)$ be one less than the number of distinct eigenvalues of this graph.

\begin{lemma}\label{lem:d_cartesian}
    We have $d(m,n,r)=m+2n$. In particular, for nonnegative integers $m_1,n_1,m_2,n_2$ with $m_h+n_h\geq 1$ for $h=1,2$,
    \[
        d(m_1+m_2,n_1+n_2,r)=d(m_1,n_1,r)+d(m_2,n_2,r).
    \]
\end{lemma}
\begin{proof}
    The eigenvalues of $K_r$ and $K_{r,r}$ are $r-1,-1$ and $r,0,-r$, respectively. Hence the distinct eigenvalues of $K_r^m\times K_{r,r}^n$ are
    \[
        m(r-1)+nr-kr\qquad (0\leq k\leq m+2n),
    \]
    so $d(m,n,r)=m+2n$, which also gives the stated identity.
\end{proof}

\begin{lemma}\label{lem:N_factors}
    For $r\geq 3$,
    \[
        N(1,0,r)=2=d(1,0,r)+1,\qquad
        N(0,1,r)=3=d(0,1,r)+1.
    \]
\end{lemma}
\begin{proof}
    For $K_r$ and $K_{r,r}$ in the orderings used in Corollary~\ref{cor:is_QP}, the second primitive idempotents are
    \[
        I_r-\frac1r J_r,\qquad
        \begin{pmatrix}
            I_r-\frac1r J_r&0\\
            0&I_r-\frac1r J_r
        \end{pmatrix},
    \]
    respectively, where $J_r$ is the all-ones matrix. Thus their distinct entries are
    \[
        \left\{-\frac1r,\frac{r-1}{r}\right\}\quad\text{and}\quad
        \left\{-\frac1r,0,\frac{r-1}{r}\right\},
    \]
    respectively. The result follows from Lemma~\ref{lem:d_cartesian}.
\end{proof}

\begin{corollary}\label{cor:not_schur_closed}
    For any integers $r\geq 3$, $m\geq 0$, and $n\geq 1$ with $m+n\geq 2$, the adjacency algebra of $K_r^m\times K_{r,r}^n$ is not closed under $\circ$.
\end{corollary}
\begin{proof}
    Fix $r\geq 3$. We prove by induction on $m+n$ that
    \[
        N(m,n,r)>d(m,n,r)+1\qquad (m\geq 0,\ n\geq 1,\ m+n\geq 2).
    \]
    For the base case $m=n=1$, by the entry formula in the proof of Proposition~\ref{prop:N_cartesian} and the preceding lemma, the distinct entries of $2r^2E_1$ for $K_r\times K_{r,r}$ are
    \[
        \{-3,-1,r-3,r-1,2r-3,3r-3\}.
    \]
    These six numbers are distinct for $r\geq 3$, so $N(1,1,r)=6>4=d(1,1,r)+1$.

    For the other base case $(m,n)=(0,2)$, the same entry formula shows that the entries of $2r^2E_1$ for $K_{r,r}^2$ are all sums of two elements of $\{-1,0,r-1\}$. The resulting values are
    \[
        \{-2,-1,0,r-2,r-1,2r-2\}.
    \]
    They are distinct for $r\geq 3$, so $N(0,2,r)=6>5=d(0,2,r)+1$.

    Now fix $k\geq 2$ and suppose the inequality holds for every pair $(m,n)$ with $m\geq 0$, $n\geq 1$, and $m+n=k$. For each such pair, Proposition~\ref{prop:N_cartesian} and Lemmas~\ref{lem:d_cartesian}--\ref{lem:N_factors} give
    \[
        \begin{aligned}
            N(m+1,n,r)&\geq N(m,n,r)+N(1,0,r)-1\\
            &>d(m,n,r)+d(1,0,r)+1\\
            &=d(m+1,n,r)+1,\\
            N(m,n+1,r)&\geq N(m,n,r)+N(0,1,r)-1\\
            &>d(m,n,r)+d(0,1,r)+1\\
            &=d(m,n+1,r)+1.
        \end{aligned}
    \]
    Every pair at level $k+1$ is covered, since decreasing $m$ if $m>0$ and decreasing $n$ otherwise gives a pair in the stated range at level $k$. This completes the induction.

    Since $E_1\in M$ is constant on each class $R_i$ defined in Section~\ref{sec:general graph}, we have $N(m,n,r)\leq\delta+1$. Thus $\delta>d(m,n,r)$, and Proposition~\ref{prop:@Ai=@Ei} shows that $M$ is not closed under $\circ$.
\end{proof}

Corollaries~\ref{cor:is_QP} and~\ref{cor:not_schur_closed} give the following theorem.

\begin{restatable}[name=\restatementtitle]{theorem}{counterexampletheorem}\label{thm:counterexample}
    For any integers $r\geq 3$, $m\geq 0$, and $n\geq 1$ with $m+n\geq 2$, the Cartesian product $K_r^m \times K_{r,r}^n$ gives a counterexample to \cite[Conjecture~20.4]{Terwilliger2023DistanceRegular}.
\end{restatable}

\begin{remark}
The complements of the counterexamples in
Theorem~\ref{thm:counterexample} are also counterexamples for the following reason: 
The counterexamples have diameter at least $3$, so their complements are
connected. A connected regular graph and its connected complement
have the same eigenspaces and adjacency algebra, and thus taking
complements preserves both the $Q$-polynomial property and the
failure of closure under entrywise multiplication for these graphs.
\end{remark}

\begin{remark}
These arguments give counterexamples of every diameter $D\geq 2$.
The choice $m=D-2$, $n=1$ in Theorem~\ref{thm:counterexample} gives diameter $D$ for $D\geq 3$, while the complement of $K_3\times K_{3,3}$ has diameter two since it is $12$-regular on $18$ vertices.
\end{remark}

In Theorem~\ref{thm:counterexample}, we construct a family of counterexamples using Cartesian products. Theorem~\ref{thm:tensor_QP} shows that we may also potentially construct counterexamples using tensor products. Meanwhile, the first remark above shows that the complement graph of each counterexample in Theorem~\ref{thm:counterexample} is again a counterexample. To see the common structure shared by all these counterexamples and potential counterexamples, we first recall the definition of an equitable partition.

\begin{definition}[Equitable partition]\label{def:equitable_partition}
A partition \(\mathcal P=\{P_1,\ldots,P_t\}\) of the vertex set of a graph is \emph{equitable} if, for every fixed \(i,j\), the number of neighbors in \(P_j\) of a vertex \(x\in P_i\) is independent of the choice of \(x\). A partition is called \emph{trivial} if it has a single part, and \emph{nontrivial} otherwise.
\end{definition}

In all these counterexamples and potential counterexamples described above, the vertex set is naturally a product of two sets, and fixing
either coordinate gives a nontrivial equitable partition. This gives us two different equitable partitions such that every part of one partition meets every
part of the other in exactly one vertex. These partitions remain
equitable after taking complements, so the same structure is
present in the complementary counterexamples.
These observations motivate the following refined conjecture.



\restatetheorem{\refinedconjecture}

We will prove Conjecture~\ref{conj:new} for abelian Cayley graphs in the next section.

\section{On the Abelian Cayley Graphs}\label{sec:new_conj}
In this section we prove Conjecture~\ref{conj:new} for abelian Cayley graphs. Let $G$ be a finite abelian group, and let $S\subseteq G$ satisfy
\[
    S=-S,\qquad
    0\notin S,\qquad
    \langle S\rangle=G.
\]

Set $\Gamma:=\Cay(G,S)$.
The graph $\Gamma$ is $Q$-polynomial over $\R$ if and only if $\Gamma$ is $Q$-polynomial over $\C$.
Indeed, since $A$ is real symmetric, its complex eigenspaces are
$V_i^{\C}=V_i\otimes_{\R}\C$.
Complexification commutes with sums and entrywise products of subspaces, and thus the $Q$-polynomiality holds over $\R$ if and only if it holds over $\C$. For this reason, we can work on eigenspaces over $\C$ instead, and we use $V_i$ to denote these complex eigenspaces in this section.

We first consider the Cayley graphs with cyclic group of prime power. Until further notice, suppose
\[
    G=\Z_{p^n},
    \qquad
    N:=p^n,
    \qquad
    \omega:=e^{2\pi i/N}.
\]
Fix the identification $\widehat G\cong G$ under which the character
indexed by $a\in G$ is given by
\[
    \chi_a(g):=\omega^{ag}.
\]
Then $\{\chi_a:a\in G\}$ forms a basis of $\C^G$, and $\chi_a$ is an
eigenvector of $\Gamma$ with eigenvalue
\[
    \theta_a:=\sum_{s\in S}\omega^{as}.
\]

Let $\lambda_0,\cdots,\lambda_d$ be an ordering of the distinct eigenvalues of $\Gamma$. Define 
\[X_i:=\{a\in G:\theta_a=\lambda_i\},\] 
then $V_i=\operatorname{span}_{\C}\{\chi_a:a\in X_i\}$ and $X_i = -X_i$. This gives the following equivalent definition of Q-polynomial.

\begin{lemma}\label{lem:QP_sumset_criterion}
    The ordering $V_0,V_1,\ldots,V_d$ is Q-polynomial if and only if, for every $0\leq i\leq d$,
    \[
        \bigsqcup_{\ell=0}^i X_\ell=\bigcup_{\ell=0}^i \ell X_1.
    \]
    Here, $X+Y:=\{x+y:x\in X,\ y\in Y\}$, and $\ell X$ denotes the $\ell$-fold sumset of $X$, with $0X:=\{0\}$.
\end{lemma}

\begin{proof}
    Since $\chi_a\circ\chi_b=\chi_{a+b}$ for all $a,b\in G$, we have
    \[
        V_1^{\circ\ell}
        =
        \operatorname{span}_{\C}
        \{\chi_a:a\in\ell X_1\}.
    \]
    Hence $\sum_{\ell=0}^iV_1^{\circ\ell}$ is spanned by the characters indexed by $\bigcup_{\ell=0}^i\ell X_1$, while $\sum_{\ell=0}^iV_\ell$ is spanned by the characters indexed by $\bigsqcup_{\ell=0}^iX_\ell$. Since the characters are linearly independent, these two subspaces are equal if and only if their indexing sets are equal. The result follows from the definition of a Q-polynomial ordering.
\end{proof}

Whenever $\Gamma$ is $Q$-polynomial, fix a $Q$-polynomial ordering of its eigenspaces. By Lemma~\ref{lem:QP_sumset_criterion}, $X_1$ generates $G$ and hence contains a unit. We choose the identification $\widehat G\cong G$ so that $1\in X_1$.

We now use the ring structure of $G=\Z_{p^n}$. For $u\in G$ and
$X\subseteq G$, write
\[
    u\cdot X:=\{ux:x\in X\}.
\]
Note that the dot distinguishes this dilation from the iterated sumset $\ell X$
defined previously.



\begin{lemma}\label{lem:Xj=uX1}
    Let $u \in G$ be a unit. If $a,b \in G$ satisfy
    $\theta_a = \theta_b$, then
    \[
        \theta_{ua} = \theta_{ub}.
    \]
    In particular, if $1 \in X_1$ and $u \in X_j$, then
    \[
        u \cdot X_1 = X_j,
        \qquad
        |X_1| = |X_j|.
    \]
\end{lemma}

\begin{proof}
    Consider the automorphism
    \[
        \sigma_u \in \operatorname{Aut}_{\mathbb Q}(\mathbb Q[\omega])
    \]
    that sends $\omega$ to $\omega^u$. For every $a \in G$,
    \[
        \theta_{ua}
        = \sum_{s \in S} \omega^{uas}
        = \sigma_u\left(\sum_{s \in S}\omega^{as}\right)
        = \sigma_u(\theta_a).
    \]
    Hence $\theta_a=\theta_b$ implies
    \[
        \theta_{ua}
        =\sigma_u(\theta_a)
        =\sigma_u(\theta_b)
        =\theta_{ub}.
    \]

    Now suppose that $1 \in X_1$ and $u \in X_j$. For every
    $a \in X_1$, we have $\theta_a=\theta_1$, and therefore
    \[
        \theta_{ua}=\theta_u.
    \]
    Thus $u\cdot X_1\subseteq X_j$. Conversely, if $b\in X_j$, then
    $\theta_b=\theta_u$. Applying the first assertion with the unit
    $u^{-1}$ gives
    \[
        \theta_{u^{-1}b}=\theta_1,
    \]
    so $u^{-1}b\in X_1$, and hence $b\in u\cdot X_1$. Therefore
    \[
        u\cdot X_1=X_j.
    \]
    Since multiplication by $u$ is a bijection of $G$, it also follows
    that $|X_1|=|X_j|$.
\end{proof}

To prepare for the proof of the next theorem, we recall Kneser's theorem
and Kemperman's structure theorem in the formulation of
\cite[Theorems A and C]{Lev06}. For a subset $C$ of an abelian group, define its period by
\[
    \pi(C):=\{g\in G:C+g=C\}.
\]
The set $C$ is called \emph{aperiodic} if $\pi(C)=\{0\}$.

\begin{theirtheorem}[Kneser]
    Let $A$ and $B$ be finite, non-empty subsets of an abelian group $G$,
satisfying $|A + B| \leq |A| + |B| - 1$. Then, letting $H := \pi(A+B)$, we have
  \[ |A+B| = |A+H|+|B+H|-|H|. \]
\end{theirtheorem}

For two finite subsets $A, B$ of an abelian group $G$ and $c \in G$, define 
\[\nu_c(A,B) := \#\{(a,b) \in A \times B : a + b = c\},\]
\[\mu(A,B) = \min\limits_{c \in A+B} \nu_c(A,B).\] We say that a pair $(A,B)$ is \emph{elementary} if at least one of the following
holds:
\begin{itemize}
\item[(a)] $\min\{|A|,|B|\}=1$;
\item[(b)] $A$ and $B$ are arithmetic progressions sharing a common
  difference $d\in G$, the order of which in $G$ is at least $|A|+|B|-1$;
\item[(c)] $A=g_1+(K_1\cup\{0\})$ and $B=g_2-(K_2\cup\{0\})$, where
  $g_1,g_2\in G$ and $K_1,K_2$ are non-empty subsets of a subgroup
  $K\le G$ such that $K=K_1\sqcup K_2\sqcup\{0\}$.
  Moreover, $c:=g_1+g_2$ is the only element of $A+B$ with $\nu_c(A,B)=1$;
\item[(d)] $A=g_1+K_1$ and $B=g_2-K_2$, where $g_1,g_2\in G$ and $K_1,K_2$
  are non-empty, aperiodic subsets of a subgroup $K\le G$ such that
  $K=K_1 \sqcup K_2$. 
  Moreover, $\mu(A,B)\ge 2$.
\end{itemize}

\begin{theirtheorem}[Kemperman structure theorem]
    Let $A$ and $B$ be finite, non-empty subsets of a non-trivial abelian group
    $G$. Then the following two conditions hold:
    \begin{itemize}
        \item[(I)] $|A + B| \leq |A| + |B| - 1$
        \item[(II)] (Kemperman's condition.) $\pi(A+B) = \{0\}$ or $\mu(A,B) = 1$.
    \end{itemize}
    if and only if there exists non-empty subsets
    $A_0\seq A,\,B_0\seq B$ and a non-zero subgroup $H\le G$ such that, writing
    $\oA:=\pgh(A)$ and $\oB:=\pgh(B)$, we have
    \begin{itemize}
        \item[(i)]   each of $A_0$ and $B_0$ is contained in an $H$-coset, and the
        pair $(A_0,B_0)$ is elementary.
        \item[(ii)]  each of $A\stm A_0$ and $B\stm B_0$ is a (possibly, empty) union
        of $H$-cosets;
        \item[(iii)] $|\oA+\oB|=|\oA|+|\oB|-1$;
        \item[(iv)]  $A_0+B_0+H$, considered as an element of $G/H$, has a unique
        representation as a sum of an element of $\oA$ and an element of $\oB$.
    \end{itemize}
\end{theirtheorem}

We will use the following special case of Kemperman Structure theorem.

\begin{corollary}\label{cor:KST}
    Let $C$ be a non-empty subset of
    $G=\Z_{p^n}$, where $p$ is prime. Suppose that
    \[
        0\in C,
        \qquad
        C=-C,
        \qquad
        \langle C\rangle=G,
    \]
    and
    \[
        \pi(2C)=\{0\},
        \qquad
        |2C|=2|C|-1.
    \]
    Then there exist a unit $c\in G$ and an integer $r\geq1$ such that
    \[
        C
        =
        \{-rc,-(r-1)c,\ldots,0,\ldots,(r-1)c,rc\}.
    \]
\end{corollary}

\begin{proof}
    Let $A = B := C$ in Kemperman structure theorem. Since \[\pi(C)\subseteq\pi(2C)=\{0\},\] the set $C$ is aperiodic. For a fixed $H$, the decomposition in Theorem~B has a unique incomplete $H$-coset: indeed, all other intersections of $C$ with $H$-cosets are complete, while $C$ cannot be a union of $H$-cosets. Since $A=B=C$, it follows that $A_0=B_0$. Denote this common set by $C_0$, and write $C_0\subseteq g+H$. We first prove that $C_0 = C$. Since $C=-C$, the unique incomplete $H$-coset is invariant under
    negation. Hence
    \[
        g+H=-g+H,
    \]
    and therefore $2\og=0$ in $G/H$. By condition (iv), \[0 = \og + \og = 0 + 0\] must be a unique representation. Hence $\og = 0$, and so $C_0 \seq H$. Moreover, for all $x \in \oC$, $0 = x + (-x)$ is another representation, and thus condition (iv) forces $\oC = \{0\}$, i.e., \[C = C_0 \seq H.\] Since $\langle C\rangle = G$, $H = G$.

    We now prove that $(C,C)$ is not of elementary type (a), (c) or (d) when $G \neq \Z_{2}$.
    \begin{itemize}
        \item [(a)] If $|C| = 1$ then $C = \{0\}$ and $\langle C \rangle \neq G$.
        \item[(c)] In this case,
        \[
            C\subseteq g_1+K=K
            \qquad\text{and}\qquad
            2C=g_1+g_2+K=K.
        \]
        Since $\langle C\rangle=G$, we obtain $K=G$, and therefore
        $2C=G$, a contradiction.
        \item [(d)] In this case \[|K| = |K_1| + |K_2| = 2|C|\] is even. Since $K$ is a subgroup of $G$ with even size, $G = \Z_{2^n}$. Note that $C = g_1+K_1 \seq g_1 + K$ and $0 \in C$, so $C \seq K$. Since $\langle C \rangle = G$, $G = K = K_1 \sqcup K_2$. Since $K_1 = -g_1 + C$ and $K_2 = g_2 - C = g_2 + C$, we have
        \[
        G = (-g_1 + C) \sqcup (g_2 + C) = C \sqcup (\gamma + C) = (\gamma + C) \sqcup (2\gamma + C)
        \]
        where $\gamma := g_1 + g_2$. Comparing the two decompositions gives
        \[
            C=2\gamma+C.
        \] 
        Since $K_1$ is aperiodic, $C = g_1 + K_1$ is aperiodic, and thus $2\gamma = 0$. Moreover, $\gamma \neq 0$ since $C \sqcup (\gamma + C)$ is a disjoint union. As a result, $\gamma = 2^{n-1}$.

        If $G = \Z_2$, then $C = \{0\}$, and thus $\mu(C,C) = 1$, contradiction. For $n \geq 2$, we can consider $x := 2^{n-2}$. Since $C = -C$, we have 
        \[x \in C \iff -x = \gamma + x \in C.\]
        On the other hand, since $C \sqcup (\gamma + C)$ is disjoint union, 
        \[x \in C \iff \gamma + x \notin C,\] 
        a contradiction.
    \end{itemize}

    Therefore, $(C,C)$ is of elementary type (b), i.e., $C$ is an arithmetic progression. Since $\langle C\rangle=G$, its common difference is a unit. The result follows from $C = -C$ and $0 \in C$.
\end{proof}

\begin{theorem}\label{thm:odd_irreducible_QP_classification}
    Suppose $\Gamma = \Cay(G,S)$ is not complete or strongly regular with some $G = \Z_{p^n}$ where $p$ is odd prime. If $\Gamma$ is Q-polynomial, then $X_i = \{\pm i\}$ for $1 \leq i \leq \frac{p^n-1}{2}$ and $X_0 = \{0\}$.
\end{theorem}

\begin{proof}
    Since $\Gamma$ is not complete or strongly regular, we have $d \geq 3$. We first prove that there exists a unit in $X_2$. Let $U$ be the collection of all units. Suppose, for the sake of contradiction, that there is no unit in $X_2$. For any $a,b \in X_1$ such that $a+b \in U$, we have
    \[a+b \in 2X_1 \seq X_0 \sqcup X_1 \sqcup X_2.\] 
    Since $X_0 = \{0\}$ and there is no unit in $X_2$, 
    \[a,b \in X_1, a+b \in U \implies a + b \in X_1.\] Since $1 \in X_1$, $2 = 1+1 \in X_1$. An induction on the
    least positive representative of $u\in U$, using
    \[
        u=(u-1)+1
        \qquad\text{or}\qquad
        u=(u-2)+2,
    \]
    shows that $U\subseteq X_1$. Since every non-unit element $x \in G$ can be expressed as the sum of two units $1$ and $x-1$, we have 
    \[X_0 \sqcup X_1 \sqcup X_2 = G,\] 
    which implies that $d \leq 2$. This shows that if $d \geq 3$ then $X_2$ contains a unit.

    By Lemma \ref{lem:Xj=uX1}, $|X_2| = |X_1|$. Now set $C:=\{0\}\cup X_1$, then
    \[
        0\in C,\qquad
        C=-C,\qquad
        \langle C\rangle=G,
    \]
    and
    \[
        |2C|
        =
        |X_0\sqcup X_1\sqcup X_2|
        =
        2|C|-1.
    \] 
    In order to use Corollary \ref{cor:KST}, it remains to check that $\pi(2C) = \{0\}$.

    Let 
    \[
        H:=\pi(2C),
        \qquad
        \oC:=(C+H)/H,
        \qquad
        t:=|\oC|,
        \qquad
        h:=|H|.
    \] 
    Then $|C + H| = th$. By Kneser's theorem, \[|2C| = 2|C + H| - |H| = 2th - h,\] and thus \[|C| = th - \frac{h - 1}{2}.\] Note that $2\oC = 2(C+H)/H = 2C/H$, so $|2\oC| = \frac{|2C|}{|H|} = 2t - 1$. 
    Since $H=\pi(2C)$, we have
    \[
        \pi(2\oC)=\{0\}.
    \]
    Moreover,
    \[
        0\in\oC,\qquad
        \oC=-\oC,\qquad
        \langle\oC\rangle=G/H,
    \]
    and
    \[
        |2\oC|=2|\oC|-1.
    \]
    Therefore, Corollary~\ref{cor:KST} implies that
    \[
        \oC
        =
        \{-r\oc,-(r-1)\oc,\ldots,0,\ldots,(r-1)\oc,r\oc\}.
    \]
    Since $2\oC \neq G/H$, $4r+1 < |G/H|$. Since $|G/H|$ is odd, $|G/H| \geq 4r+3$. Note that $\oc$ is a unit in $G/H$ since $\oC$ generates $G/H$, so both $(2r+1)\oc$ and $(2r+2)\oc$ are not in $2\oC$, and at least one of them is a unit. Denote the unit as $\ou$. Let $u \in U$ be a unit in $\pgh^{-1}(\ou)$, then $u \in X_j$ for some $j \geq 3$. By Lemma \ref{lem:Xj=uX1}, $X_j = u \cdot X_1$.

    Suppose, for the sake of contradiction that $|H| > 1$. Since 
    \[|C \cap H| = |H| - |H\stm C| \geq |H| - |(C+H) \stm C| = |H| - |C+H| + |C| = \frac{h+1}{2} > 1,\]
    there exists $0 \neq x \in X_1 \cap H$. Since $H$ is invariant under a multiplication of a unit, $ux \in H$. Moreover, as proved previously, $X_j = u \cdot X_1$, which implies that $ux \in X_j$, so $X_j \cap H \neq \emptyset$. However, $0 \in 2C$ implies that $H = \pi(2C) \seq 2C$. Since $2C \cap X_j = \emptyset$ for $j \geq 3$, $X_j \cap H = \emptyset$, a contradiction. Therefore, $|H| = 1$ and it follows that $H = \{0\}$.

    By Corollary \ref{cor:KST}, $C$ is of the form $C = \{-rc, -(r-1)c, \dots, 0, \dots, (r-1)c, rc\}$. Since $\langle C \rangle = G$, $c \in U \cap X_1$. By Lemma \ref{lem:Xj=uX1}, $c \cdot X_1 = X_1$, so
    \[X_1 = c^{-1} \cdot X_1 = \{\pm 1, \pm 2, \dots, \pm r\}.\]
    We now prove that $X_i = \{\pm i\}$ for all $i$. Suppose, for the sake of contradiction, that $r \geq 2$. Consider
    \[x = \begin{cases}
        r, &p \nmid r,\\
        r-1, &p \mid r.
    \end{cases}\]
    Then $x \in U \cap X_1$. By Lemma \ref{lem:Xj=uX1} and since $2 \in X_1$, $2x \in X_1$. However, since $\pi(2C) = \{0\}$, $2C \neq G$, i.e., $4r+1 < p^n$. Moreover, since $p \neq 2$, $r=2$ implies $x=r$. As a consequence, $r < 2x \leq 2r < \frac{p^n}{2}$, and thus $2x \notin X_1$, contradiction. Therefore, \[X_1 = \{\pm 1\}\] and thus by Lemma~\ref{lem:QP_sumset_criterion} \[X_i = (\bigcup_{\ell = 0}^i \ell X_1) \stm \bigsqcup_{\ell = 0}^{i-1} X_\ell = \{\pm i\}, 1 
    \leq i \leq \frac{N - 1}{2}\]
    Since these already covers the whole $G$, these are all of the $X_i$.
    
\end{proof}



\begin{theorem}\label{thm:even_irreducible_QP_classification}
    Suppose $\Gamma = \Cay(G,S)$ is not complete or strongly regular with some $G = \Z_{2^n}$. If $\Gamma$ is Q-polynomial, then $X_i = \{\pm i\}$ for $1 \leq i < 2^{n-1}$, $X_0 = \{0\}$, and $X_{2^{n-1}} = \{2^{n-1}\}$.
\end{theorem}
\begin{proof}
    Since $\Gamma$ is not complete or strongly regular, we have $d \geq 3$. Let $U$ be the collection of all units in $G$. We first show that $X_2 \cap U = \emptyset$ by following the proof of Theorem~\ref{thm:odd_irreducible_QP_classification}. 
    
    Suppose, for a contradiction, that $X_2 \cap U \neq \emptyset$. By Lemma~\ref{lem:Xj=uX1}, we have $|X_2|=|X_1|$. Setting \(C:=\{0\}\cup X_1\), we obtain, as in the proof of Theorem~\ref{thm:odd_irreducible_QP_classification},
    \[
    |2C|=2|C|-1.
    \]
    By definition, $2C$ is a disjoint union of $\pi(2C)$-cosets, so $|\pi(2C)| \mid |2C|$, and thus $|\pi(2C)|$ is odd. Moreover, $\pi(2C)$ is a subgroup of $G$, so $|\pi(2C)| \mid 2^n$. Therefore, $|\pi(2C)| = 1$ and $\pi(2C) = \{0\}$. We may therefore apply Corollary~\ref{cor:KST} and follow the normalization argument in the proof of Theorem~\ref{thm:odd_irreducible_QP_classification} to obtain
    \[
    X_1=\{\pm1,\ldots,\pm r\},
    \qquad 4r+1<2^n.
    \]
    We now prove that $r$ = 1. The same argument as in the proof of Theorem 4.4 excludes $r\geq3$. For $r=2$, the sumset criterion for $Q$-polynomial (Lemma~\ref{lem:QP_sumset_criterion}) gives that \[X_2 = \{\pm 3, \pm 4\}.\] By Lemma \ref{lem:Xj=uX1}, \[X_2 = 3 \cdot X_1 = \{\pm 3, \pm 6\}.\] Since $9 = 4r+1 < 2^n$, we have $2^n \geq 16$, which shows that $6 \notin \{\pm 3, \pm 4\}$, a contradiction. Therefore, $r=1$ and $X_2 = \{\pm 2\}$, which indicates that $X_2 \cap U = \emptyset$, a contradiction. Therefore, $X_2 \cap U = \emptyset$. 

    We next show that $X_1 \seq U$. Otherwise, fix any $a \in X_1 \setminus U$ and put $K := X_1 \cap U$. For any $k \in K$, $k + a \in U$. Since $k,a \in X_1$, by Lemma~\ref{lem:QP_sumset_criterion}
    \[k+a \in 2X_1 \seq X_0 \sqcup X_1 \sqcup X_2.\]
    Since $U$ intersects with neither $X_0 = \{0\}$ nor $X_2$, $k+a \in X_1$, and thus $k+a \in K$. Therefore, $K+a = K$. In particular, $-1 \in K$, so for every \(u\in U\) we have \(-1+ua\in K\), and thus
    \[
    ua=1+(-1+ua)\in2X_1.
    \]
    Moreover, $\theta_a=\theta_1$, so Lemma~\ref{lem:Xj=uX1} gives $\theta_{ua}=\theta_u$. It follows that $u\in X_1$. Since $u$ is arbitrary, $U \seq X_1$. Every even element of $G$ is a sum of two units, so $G = X_0 \sqcup X_1 \sqcup X_2$ and thus $d \leq 2$, contrary to the fact that $\Gamma$ is not complete or strongly regular. Thus $X_1 \seq U$.

    For any $a \in X_1$, Lemma~\ref{lem:Xj=uX1} now gives that $aX_1 = X_1$, which implies that $X_1$ is a multiplication subgroup of $U$. Since $X_1$ is a multiplicative subgroup of \(U\) containing $-1$, the explicit classification of these subgroups in \cite[p.~3, Eq.~(1)]{Kovacs2005} gives
    \[
    X_1=(1+2^r\cdot G)\cup(-1+2^r\cdot G)
    \]
    for some \(2\le r\le n\). By Lemma~\ref{lem:QP_sumset_criterion}, we have
    \[X_i = \begin{cases}
        \{0\}, \quad &i = 0,\\
        (2 + 2^r \cdot G) \cup (2^r \cdot G) \cup (-2 + 2^r \cdot G) \stm \{0\}, \quad &i = 2,\\
        (i + 2^r \cdot G) \cup (-i + 2^r \cdot G), \quad &\text{otherwise}.
    \end{cases}\]
    Therefore, $X_1$, $X_2 \cup \{0\}$, and every $X_i$ with $i \geq 3$ are unions of $H$-cosets, where $H := 2^r \cdot G$.

    We now prove that $r = n$. Suppose, for the contradiction, that $r < n$. 
    Let $s \in S \cap U$, which exists since $\langle S \rangle =G$. For any $H$-coset $b + H$, we have
    \[
    \sum_{a\in b+H}\omega^{-as}
    =
    \omega^{-bs}
    \sum_{j=0}^{2^{n-r}-1}
    \left(\omega^{-2^rs}\right)^j
    =0.
    \]
    The last equality follows from the fact that \(\omega^{-2^rs}\) is a nontrivial \(2^{n-r}\)-th root of unity. Consequently,
    \[
    \sum_{a\in X_i}\omega^{-as}
    =
    \begin{cases}
    1,&i=0,\\
    -1,&i=2,\\
    0,&\text{otherwise}.
    \end{cases}
    \]
    By definition, for any $a \in G$,
    \[
    \theta_a = \sum_{g \in G} \mathbbm{1}_S(g) \omega^{ag}.
    \]
    By Fourier inverse formula,
    \[
    1 = \mathbbm{1}_S(s) = \frac{1}{2^n} \sum_{a \in G} \theta_a \omega^{-as} = \frac{1}{2^n} \sum_{i=0}^{d} \lambda_i \sum_{a\in X_i} \omega^{-as} = \frac{\lambda_0 - \lambda_2}{2^n},
    \]
    and hence $\lambda_2 = \lambda_0 - 2^n$. 

    Let \(\overline{S}:=G\setminus(S\cup\{0\})\), so that the complement of \(\Gamma\) is \(\overline{\Gamma}=\operatorname{Cay}(G,\overline{S})\). The eigenvalue of \(\overline{\Gamma}\) corresponding to the nontrivial character \(\chi_2\) is therefore
    \[
    \overline{\theta}_2
    =\sum_{\overline{s}\in\overline{S}}\omega^{2\overline{s}}
    =-1-\theta_2
    =-1-\lambda_2
    =2^n-|S|-1
    =|\overline{S}|.
    \]
    Therefore, $\omega^{2\overline{s}} = 1$ for any $\overline{s} \in \overline{S}$, and hence $\overline{S} = \emptyset$ or $\{2^{n-1}\}$. It follows that \(\Gamma\) is complete in the former case and strongly regular in the latter case, contradicting the hypothesis. Therefore, $r = n$ and thus
    \[
    X_i = \begin{cases}
        \{0\}, \quad &i=0,\\
        \{2^{n-1}\}, \quad &i = 2^{n-1},\\
        \{\pm i\}, \quad &\text{otherwise}.
    \end{cases}
    \]

\end{proof}

Theorem~\ref{thm:odd_irreducible_QP_classification} and Theorem~\ref{thm:even_irreducible_QP_classification} directly yield the following corollary.

\begin{corollary}\label{cor:irreducible_QP_SC}
    If $\Gamma = \Cay(G,S)$ is Q-polynomial for some $G = \Z_{p^n}$, then the adjacency algebra $\R[A]$ is closed under entrywise multiplication $\circ$.
\end{corollary}
\begin{proof}
    The complete and strongly regular cases follow from Section 2.3. In the remaining cases, put \(N=p^n\). By Theorems 4.4 and 4.5, the character classes of \(\Gamma\) are the singleton \(\{0\}\), the two-element sets \(\{i,-i\}\) for \(1\le i<N/2\), and, when \(N\) is even, the additional singleton \(\{N/2\}\).

    For the cycle \(C_N\), the eigenvalue on \(\chi_i\) is \(2\cos(2\pi i/N)\), which gives exactly the same character classes. Thus \(\Gamma\) and \(C_N\) have the same eigenspaces and hence the same adjacency algebra. Since \(C_N\) is distance-regular, this algebra is closed under entrywise multiplication.
\end{proof}

Now we consider the general Abelian Cayley graphs. Suppose $G$ is a general abelian group.
\begin{lemma}\label{lem:reducible_case_satisfy}
    If there exists a nontrivial decomposition $G = G_1 \oplus G_2$, then there exist two nontrivial equitable partitions $\{\cP_i\}, \{\cQ_j\}$ of the vertices, such that $|\cP_i \cap \cQ_j| = 1$ for all $i,j$.
\end{lemma}

\begin{proof}
    Let $G_1$ and $G_2$ act on $G=G_1\oplus G_2$ by translation in the first and second coordinates, respectively. For $i\in G_2$ and $j\in G_1$, let
    \[
        \begin{aligned}
            \cP_i&:=\operatorname{Orb}_{G_1}(0,i)=\{(g,i):g\in G_1\},\\
            \cQ_j&:=\operatorname{Orb}_{G_2}(j,0)=\{(j,h):h\in G_2\}.
        \end{aligned}
    \]
    These actions induce the partitions
    $\{\cP_i\}_{i\in G_2}$ and $\{\cQ_j\}_{j\in G_1}$. Since $|G_1|,|G_2|>1$, each partition has at least two parts, each of size at least two; hence both partitions are nontrivial. Moreover,
    \[
        \cP_i\cap \cQ_j=\{(j,i)\},
        \qquad
        |\cP_i\cap \cQ_j|=1.
    \]
    for all $i\in G_2$ and $j\in G_1$.

    We show that $\{\cP_i\}_{i\in G_2}$ is equitable. Fix $i,k\in G_2$ and $x,x'\in \cP_i$. Since $x$ and $x'$ lie in the same $G_1$-orbit, there exists $a\in G_1$ such that $x'=x+(a,0)$. Translation by $(a,0)$ is an automorphism of $\Gamma$ and preserves $\cP_k$ setwise. It therefore induces a bijection from $\Gamma(x)\cap \cP_k$ to $\Gamma(x')\cap \cP_k$. Thus $|\Gamma(x)\cap \cP_k|$ is independent of the choice of $x\in \cP_i$, so $\{\cP_i\}_{i\in G_2}$ is equitable. The same argument shows that $\{\cQ_j\}_{j\in G_1}$ is equitable.
\end{proof}

By the structure theorem for finite abelian groups, every nontrivial finite abelian group is either cyclic of prime-power order or admits a nontrivial direct sum decomposition. Corollary~\ref{cor:irreducible_QP_SC} and Lemma~\ref{lem:reducible_case_satisfy} therefore establish Conjecture~\ref{conj:new} in these two cases.

\restatetheorem{\abeliancayleytheorem}

\section*{Acknowledgments}
I am grateful to my advisor, Paul Terwilliger, for his guidance,
helpful discussions, and valuable suggestions throughout this work.

Generative AI tools were used to assist with mathematical
exploration, proof development, and language editing.
The author takes full responsibility for the mathematical
arguments and the final content of this paper.

\bibliographystyle{plain}
\enlargethispage{4\baselineskip}
\bibliography{reference}

\end{document}